\documentclass[11pt]{amsart}
\usepackage{amsmath,amssymb}
\newtheorem{theorem}{Theorem}[section]
\newtheorem{proposition}[theorem]{Proposition}
\newtheorem{corollary}[theorem]{Corollary}

\begin{document}

\title[Mean curvature flow in Riemannian manifolds]{An inscribed radius estimate for mean curvature flow in Riemannian manifolds}
\author{Simon Brendle}
\address{Department of Mathematics \\ Stanford University \\ Stanford, CA 94305}
\thanks{The author was supported in part by the National Science Foundation under grant DMS-1201924.}
\begin{abstract}
We consider a family of embedded, mean convex hypersurfaces in a Riemannian manifold which evolve by the mean curvature flow. We show that, given any number $T>0$ and any $\delta>0$, we can find a constant $C_0$ with the following property: if $t \in [0,T)$ and $p$ is a point on $M_t$ where the curvature is greater than $C_0$, then the inscribed radius is at least $\frac{1}{(1+\delta) \, H}$ at the point $p$. The constant $C_0$ depends only on $\delta$, $T$, and the initial data.
\end{abstract}
\maketitle 

\section{Introduction} 

Our goal in this paper is to generalize in \cite{Brendle} to solutions of the mean curvature flow in a Riemannian manifold. Let $N$ be a Riemannian manifold of dimension $n+1$, and let $F: M \times [0,T) \to N$ be a family of closed, embedded, mean convex hypersurfaces in $N$ which evolve by mean curvature flow. As in \cite{Brendle}, we define a function $\mu$ by 
\[\mu(x,t) = \sup_{y \in M, \, 0 < d(F(x,t),F(y,t)) \leq \frac{1}{2} \, \text{\rm inj}(N)} \Big ( -\frac{2 \, \langle \exp_{F(x,t)}^{-1}(F(y,t)),\nu(x,t) \rangle}{d(F(x,t),F(y,t))^2} \Big ).\] 
Note that $\lambda_1 \leq \hdots \leq \lambda_n \leq \mu$, where the $\lambda_i$ are the principal curvatures. The reciprocal of $\mu(x,t)$ can be interpreted as the inscribed radius of $M_t$ at the point $x$.

\begin{theorem} 
\label{inscribed.radius}
Let $\delta > 0$ and $T > 0$ be given positive numbers. Then the function $\mu$ satisfies an estimate of the form 
\[\mu \leq (1+\delta) \, H + C(N,M_0,\delta,T)\] 
for all $t \in [0,T)$ and all points on $M_t$. 
\end{theorem}

In the special case that $N$ is the Euclidean space $\mathbb{R}^{n+1}$, it follows from general results of Brian White that the ratio $\frac{\mu}{H}$ is uniformly bounded from above (cf. \cite{White1}, \cite{White2}, \cite{White3}). Later, Andrews \cite{Andrews} gave a direct proof of that fact using the maximum principle. In a recent paper \cite{Brendle}, we showed that, for any mean convex solution to the mean curvature flow in Euclidean space, we have an estimate of the form $\mu \leq (1+\delta) \, H + C$, where $C$ is a positive constant that depends only on $\delta$ and the initial hypersurface $M_0$. Theorem \ref{inscribed.radius} generalizes this result to Riemannian manifolds.

We next define 
\[\rho(x,t) = \max \bigg \{ \sup_{y \in M, \, 0 < d(F(x,t),F(y,t)) \leq \frac{1}{2} \, \text{\rm inj}(N)} \frac{2 \, \langle \exp_{F(x,t)}^{-1}(F(y,t)),\nu(x,t) \rangle}{d(F(x,t),F(y,t))^2},0 \bigg \}.\]
Note that $-\rho \leq \lambda_1 \leq \hdots \leq \lambda_n$. The reciprocal of $\rho(x,t)$ has a geometric interpretation as the outer radius of $M_t$ at the point $x$.

\begin{theorem} 
\label{outer.radius}
Let $\delta > 0$ and $T > 0$ be given positive numbers. Then the function $\rho$ satisfies an estimate of the form 
\[\rho \leq \delta \, H + C(N,M_0,\delta,T)\] 
for all $t \in [0,T)$ and all points on $M_t$. 
\end{theorem}

We note that Theorem \ref{outer.radius} is a refinement of the convexity estimate of Huisken and Sinestrari \cite{Huisken-Sinestrari1}, \cite{Huisken-Sinestrari2}; see also \cite{White1}, \cite{White2}, \cite{White3}. 

\section{Evolution of the inscribed radius under mean curvature flow}

\label{calc}

Given any point $q \in N$, we define a function $\psi_q: N \to \mathbb{R}$ by $\psi_q(p) = \frac{1}{2} \, d(p,q)^2$, where $d(p,q)$ denotes the Riemannian distance in $N$. Moreover, let us put $\Xi_{q,p} := (\text{\rm Hess} \, \psi_q)_p - g$. Clearly, $\Xi_{q,p}$ is a symmetric bilinear form on $T_p N$, and we have $|\Xi_{q,p}| \leq O(d(p,q)^2)$.

\begin{proposition}
\label{evol}
Consider a point $(\bar{x},\bar{t}) \in M \times [0,T)$ such that $\lambda_n(\bar{x},\bar{t}) < \mu(\bar{x},\bar{t})$ and $\mu(\bar{x},\bar{t}) \geq 8 \, \text{\rm inj}(N)^{-1}$. We further assume that $U$ is an open neighborhood of $\bar{x}$ and $\Phi: U \times (\bar{t}-\alpha,\bar{t}] \to \mathbb{R}$ is a smooth function such that $\Phi(\bar{x},\bar{t}) = \mu(\bar{x},\bar{t})$ and $\Phi(x,t) \geq \mu(x,t)$ for all points $(x,t) \in U \times (\bar{t}-\alpha,\bar{t}]$. Then 
\[\frac{\partial \Phi}{\partial t} - \Delta \Phi - |A|^2 \, \Phi + \sum_{i=1}^n \frac{1}{\Phi-\lambda_i} \, (D_i \Phi)^2 \leq C \, H + C \, \Phi + C \, \sum_{i=1}^n \frac{1}{\Phi-\lambda_i}\] 
at the point $(\bar{x},\bar{t})$. Here, $C$ is a positive constant that depends only on the ambient manifold $N$ and the initial hypersurface $M_0$.
\end{proposition}

\textbf{Proof.} 
Let us define a function $Z: M \times M \times [0,T) \to \mathbb{R}$ by 
\begin{align*} 
Z(x,y,t) 
&= \Phi(x,t) \, \psi_{F(y,t)}(F(x,t)) - \big \langle \nabla \psi_{F(y,t)} \big |_{F(x,t)},\nu(x,t) \big \rangle \\ 
&= \frac{1}{2} \, \Phi(x,t) \, d(F(x,t),F(y,t))^2 + \big \langle \exp_{F(x,t)}^{-1}(F(y,t)),\nu(x,t) \big \rangle. 
\end{align*} 
By assumption, we have $Z(x,y,t) \geq 0$ whenever $x \in U$, $t \in (\bar{t}-\alpha,\bar{t}]$, and $d(F(x,t),F(y,t)) \leq \frac{1}{2} \, \text{\rm inj}(N)$. Moreover, we can find a point $\bar{y} \in M$ such that $0 < d(F(\bar{x},\bar{t}),F(\bar{y},\bar{t})) \leq \frac{1}{2} \, \text{\rm inj}(N)$ and $Z(\bar{x},\bar{y},\bar{t}) = 0$. It is clear that $\Phi(\bar{x},\bar{t}) \, d(F(\bar{x},\bar{t}),F(\bar{y},\bar{t})) \leq 2$, so $d(F(\bar{x},\bar{t}),F(\bar{y},\bar{t})) \leq \frac{1}{4} \, \text{\rm inj}(N)$. This implies 
\begin{align*}
0 = \frac{\partial Z}{\partial x_i}(\bar{x},\bar{y},\bar{t}) 
&= \frac{1}{2} \, \frac{\partial \Phi}{\partial x_i}(\bar{x},\bar{t}) \, d(F(\bar{x},\bar{t}),F(\bar{y},\bar{t}))^2 \\ 
&- \Phi(\bar{x},\bar{t}) \, \Big \langle \exp_{F(\bar{x},\bar{t})}^{-1}(F(\bar{y},\bar{t})),\frac{\partial F}{\partial x_i}(\bar{x},\bar{t}) \Big \rangle \\ 
&+ h_i^j(\bar{x},\bar{t}) \, \Big \langle \exp_{F(\bar{x},\bar{t})}^{-1}( F(\bar{y},\bar{t})),\frac{\partial F}{\partial x_j}(\bar{x},\bar{t}) \Big \rangle \\ 
&- \Xi_{F(\bar{y},\bar{t}),F(\bar{x},\bar{t})} \Big ( \frac{\partial F}{\partial x_i}(\bar{x},\bar{t}),\nu(\bar{x},\bar{t}) \Big ). 
\end{align*} 
Rearranging terms gives 
\begin{align*} 
&\Big \langle \exp_{F(\bar{x},\bar{t})}^{-1}( F(\bar{y},\bar{t})),\frac{\partial F}{\partial x_i}(\bar{x},\bar{t}) \Big \rangle \\ 
&= \frac{1}{2} \, \frac{1}{\Phi(\bar{x},\bar{t}) - \lambda_i(\bar{x},\bar{t})} \, \Big ( \frac{\partial \Phi}{\partial x_i}(\bar{x},\bar{t}) + O(1) \Big ) \, d(F(\bar{x},\bar{t}),F(\bar{y},\bar{t}))^2. 
\end{align*} 
We now differentiate one more time. Using the Codazzi equations, we obtain 
\begin{align*}
&\sum_{i=1}^n \frac{\partial^2 Z}{\partial x_i^2}(\bar{x},\bar{y},\bar{t}) \\ 
&= \Delta \Phi(\bar{x},\bar{t}) \, \psi_{F(\bar{y},\bar{t})}(F(\bar{x},\bar{t})) \\ 
&- 2 \, \frac{\partial \Phi}{\partial x_i}(\bar{x},\bar{t}) \, \Big \langle \exp_{F(\bar{x},\bar{t})}^{-1}( F(\bar{y},\bar{t})),\frac{\partial F}{\partial x_i}(\bar{x},\bar{t}) \Big \rangle \\ 
&+ \frac{\partial H}{\partial x_i}(\bar{x},\bar{t}) \, \Big \langle \exp_{F(\bar{x},\bar{t})}^{-1}( F(\bar{y},\bar{t})),\frac{\partial F}{\partial x_i}(\bar{x},\bar{t}) \Big \rangle \\ 
&+ H(\bar{x},\bar{t}) \, \Phi(\bar{x},\bar{t}) \, \langle \exp_{F(\bar{x},\bar{t})}^{-1}( F(\bar{y},\bar{t})),\nu(\bar{x},\bar{t}) \rangle \\
&- |A(\bar{x},\bar{t})|^2 \, \langle \exp_{F(\bar{x},\bar{t})}^{-1}( F(\bar{y},\bar{t})),\nu(\bar{x},\bar{t}) \rangle \\  
&+ n \, \Phi(\bar{x},\bar{t}) - H(\bar{x},\bar{t}) \\ 
&+ O \big ( d(F(\bar{x},\bar{t}),F(\bar{y},\bar{t})) + H(\bar{x},\bar{t}) \, d(F(\bar{x},\bar{t}),F(\bar{y},\bar{t}))^2 \big ), 
\end{align*} 
hence 
\begin{align*}
&\sum_{i=1}^n \frac{\partial^2 Z}{\partial x_i^2}(\bar{x},\bar{y},\bar{t}) \\ 
&\leq \frac{1}{2} \, \bigg ( \Delta \Phi(\bar{x},\bar{t}) + |A(\bar{x},\bar{t})|^2 \, \Phi(\bar{x},\bar{t}) \\ 
&\hspace{20mm} - \sum_{i=1}^n \frac{2}{\Phi(\bar{x},\bar{t})-\lambda_i(\bar{x},\bar{t})} \, \Big ( \frac{\partial \Phi}{\partial x_i}(\bar{x},\bar{t}) \Big )^2 \bigg ) \, d(F(\bar{x},\bar{t}),F(\bar{y},\bar{t}))^2 \\ 
&+ \frac{\partial H}{\partial x_i}(\bar{x},\bar{t}) \, \Big \langle \exp_{F(\bar{x},\bar{t})}^{-1}(F(\bar{y},\bar{t})),\frac{\partial F}{\partial x_i}(\bar{x},\bar{t}) \Big \rangle \\ 
&+ H(\bar{x},\bar{t}) \, \Phi(\bar{x},\bar{t}) \, \big \langle \exp_{F(\bar{x},\bar{t})}^{-1}(F(\bar{y},\bar{t})),\nu(\bar{x},\bar{t}) \big \rangle \\
&+ n \, \Phi(\bar{x},\bar{t}) - H(\bar{x},\bar{t}) \\ 
&+ O \big ( d(F(\bar{x},\bar{t}),F(\bar{y},\bar{t})) + H(\bar{x},\bar{t}) \, d(F(\bar{x},\bar{t}),F(\bar{y},\bar{t}))^2 \big ) \\ 
&+ O \bigg ( \sum_{i=1}^n \frac{1}{\Phi(\bar{x},\bar{t})-\lambda_i(\bar{x},\bar{t})} \, \Big | \frac{\partial \Phi}{\partial x_i}(\bar{x},\bar{t}) \Big | \, d(F(\bar{x},\bar{t}),F(\bar{y},\bar{t}))^2 \bigg ). 
\end{align*} 
We next compute 
\begin{align*} 
0 &= \frac{\partial Z}{\partial y_i}(\bar{x},\bar{y},\bar{t}) \\ 
&= \Big \langle (D\exp_{F(\bar{x},\bar{t})}^{-1})_{F(\bar{y},\bar{t})} \Big ( \frac{\partial F}{\partial y_i}(\bar{y},\bar{t}) \Big ),\nu(\bar{x},\bar{t}) + \Phi(\bar{x},\bar{t}) \, \exp_{F(\bar{x},\bar{t})}^{-1}(F(\bar{y},\bar{t})) \Big \rangle. 
\end{align*} 
This implies 
\begin{align*} 
&\sum_{i=1}^n \frac{\partial^2 Z}{\partial y_i^2}(\bar{x},\bar{y},\bar{t}) \\ 
&= n \, \Phi(\bar{x},\bar{t}) - H(\bar{y},\bar{t}) \, \big \langle (D\exp_{F(\bar{x},\bar{t})}^{-1})_{F(\bar{y},\bar{t})}(\nu(\bar{y},\bar{t})),\nu(\bar{x},\bar{t}) + \Phi(\bar{x},\bar{t}) \, \exp_{F(\bar{x},\bar{t})}^{-1}(F(\bar{y},\bar{t})) \big \rangle \\ 
&+ O \big ( d(F(\bar{x},\bar{t}),F(\bar{y},\bar{t})) \big ). 
\end{align*} 
Note that the vector $\nu(\bar{x},\bar{t}) + \Phi(\bar{x},\bar{t}) \, \exp_{F(\bar{x},\bar{t})}^{-1}(F(\bar{y},\bar{t}))$ has unit length. From this, we deduce that 
\begin{align*} 
(D\exp_{F(\bar{x},\bar{t})}^{-1})_{F(\bar{y},\bar{t})} (\nu(\bar{y},\bar{t})) 
&= \nu(\bar{x},\bar{t}) + \Phi(\bar{x},\bar{t}) \, \exp_{F(\bar{x},\bar{t})}^{-1}(F(\bar{y},\bar{t})) \\ 
&+ O \big ( d(F(\bar{x},\bar{t}),F(\bar{y},\bar{t}))^2 \big ). 
\end{align*} 
Moreover, we can arrange that 
\begin{align*} 
&(D\exp_{F(\bar{x},\bar{t})}^{-1})_{F(\bar{y},\bar{t})} \Big ( \frac{\partial F}{\partial y_i}(\bar{y},\bar{t}) \Big ) \\ 
&= \frac{\partial F}{\partial x_i}(\bar{x},\bar{t}) - 2 \, \frac{\langle \exp_{F(\bar{x},\bar{t})}^{-1}(F(\bar{y},\bar{t})),\frac{\partial F}{\partial x_i}(\bar{x},\bar{t}) \rangle}{d(F(\bar{x},\bar{t}),F(\bar{y},\bar{t}))^2} \, \frac{\exp_{F(\bar{x},\bar{t})}^{-1}(F(\bar{y},\bar{t}))}{d(F(\bar{x},\bar{t}),F(\bar{y},\bar{t}))^2} \\ 
&+ O \big ( d(F(\bar{x},\bar{t}),F(\bar{y},\bar{t}))^2 \big ). 
\end{align*} 
This gives  
\begin{align*} 
\frac{\partial^2 Z}{\partial x_i \, \partial y_i}(\bar{x},\bar{y},\bar{t}) 
&= \frac{\partial \Phi}{\partial x_i}(\bar{x},\bar{t}) \, \Big \langle (D\exp_{F(\bar{x},\bar{t})}^{-1})_{F(\bar{y},\bar{t})} \Big ( \frac{\partial F}{\partial y_i}(\bar{y},\bar{t}) \Big ),\exp_{F(\bar{x},\bar{t})}^{-1}(F(\bar{y},\bar{t})) \Big \rangle \\ 
&- (\Phi(\bar{x},\bar{t})-\lambda_i(\bar{x},\bar{t})) \, \Big \langle (D\exp_{F(\bar{x},\bar{t})}^{-1})_{F(\bar{y},\bar{t})} \Big ( \frac{\partial F}{\partial y_i}(\bar{y},\bar{t}) \Big ),\frac{\partial F}{\partial x_i}(\bar{x},\bar{t}) \Big \rangle \\ 
&+ O \big ( d(F(\bar{x},\bar{t}),F(\bar{y},\bar{t})) \big ) \\ 
&= -\frac{\partial \Phi}{\partial x_i}(\bar{x},\bar{t}) \, \Big \langle \exp_{F(\bar{x},\bar{t})}^{-1}(F(\bar{y},\bar{t})),\frac{\partial F}{\partial x_i}(\bar{x},\bar{t}) \Big \rangle \\ 
&- (\Phi(\bar{x},\bar{t})-\lambda_i(\bar{x},\bar{t})) \, \Big ( 1 - 2 \,\frac{\langle \exp_{F(\bar{x},\bar{t})}^{-1}(F(\bar{y},\bar{t})),\frac{\partial F}{\partial x_i}(\bar{x},\bar{t}) \rangle^2}{d(F(\bar{x},\bar{t}),F(\bar{y},\bar{t}))^2} \Big ) \\ 
&+ O \big ( d(F(\bar{x},\bar{t}),F(\bar{y},\bar{t})) + H(\bar{x},\bar{t}) \, d(F(\bar{x},\bar{t}),F(\bar{y},\bar{t}))^2 \big ) \\ 
&= -(\Phi(\bar{x},\bar{t})-\lambda_i(\bar{x},\bar{t})) \\ 
&+ O \big ( d(F(\bar{x},\bar{t}),F(\bar{y},\bar{t})) + H(\bar{x},\bar{t}) \, d(F(\bar{x},\bar{t}),F(\bar{y},\bar{t}))^2 \big )
\end{align*} 
for each $i$. Summation over $i$ gives 
\begin{align*} 
\sum_{i=1}^n \frac{\partial^2 Z}{\partial x_i \, \partial y_i}(\bar{x},\bar{y},\bar{t}) 
&= -n \, \Phi(\bar{x},\bar{t}) + H(\bar{x},\bar{t}) \\ 
&+ O \big ( d(F(\bar{x},\bar{t}),F(\bar{y},\bar{t})) + H(\bar{x},\bar{t}) \, d(F(\bar{x},\bar{t}),F(\bar{y},\bar{t}))^2 \big ). 
\end{align*} 
Thus, we conclude that 
\begin{align*}
&\sum_{i=1}^n \Big ( \frac{\partial^2 Z}{\partial x_i^2}(\bar{x},\bar{y},\bar{t}) + 2 \, \frac{\partial^2 Z}{\partial x_i \, \partial y_i}(\bar{x},\bar{y},\bar{t}) + \frac{\partial^2 Z}{\partial y_i^2}(\bar{x},\bar{y},\bar{t}) \Big ) \\ 
&\leq \frac{1}{2} \, \bigg ( \Delta \Phi(\bar{x},\bar{t}) + |A(\bar{x},\bar{t})|^2 \, \Phi(\bar{x},\bar{t}) \\ 
&\hspace{20mm} - \sum_{i=1}^n \frac{2}{\Phi(\bar{x},\bar{t})-\lambda_i(\bar{x},\bar{t})} \, \Big ( \frac{\partial \Phi}{\partial x_i}(\bar{x},\bar{t}) \Big )^2 \bigg ) \, d(F(\bar{x},\bar{t}),F(\bar{y},\bar{t}))^2 \\ 
&+ \frac{\partial H}{\partial x_i}(\bar{x},\bar{t}) \, \Big \langle \exp_{F(\bar{x},\bar{t})}^{-1}(F(\bar{y},\bar{t})),\frac{\partial F}{\partial x_i}(\bar{x},\bar{t}) \Big \rangle \\ 
&+ H(\bar{x},\bar{t}) + H(\bar{x},\bar{t}) \, \Phi(\bar{x},\bar{t}) \, \big \langle \exp_{F(\bar{x},\bar{t})}^{-1}(F(\bar{y},\bar{t})),\nu(\bar{x},\bar{t}) \big \rangle \\
&- H(\bar{y},\bar{t}) \, \big \langle (D\exp_{F(\bar{x},\bar{t})}^{-1})_{F(\bar{y},\bar{t})}(\nu(\bar{y},\bar{t})),\nu(\bar{x},\bar{t}) + \Phi(\bar{x},\bar{t}) \, \exp_{F(\bar{x},\bar{t})}^{-1}(F(\bar{y},\bar{t})) \big \rangle \\ 
&+ O \big ( d(F(\bar{x},\bar{t}),F(\bar{y},\bar{t})) + H(\bar{x},\bar{t}) \, d(F(\bar{x},\bar{t}),F(\bar{y},\bar{t}))^2 \big ) \\ 
&+ O \bigg ( \sum_{i=1}^n \frac{1}{\Phi(\bar{x},\bar{t})-\lambda_i(\bar{x},\bar{t})} \, \Big | \frac{\partial \Phi}{\partial x_i}(\bar{x},\bar{t}) \Big | \, d(F(\bar{x},\bar{t}),F(\bar{y},\bar{t}))^2 \bigg ). 
\end{align*} 
On the other hand, we have 
\begin{align*} 
\frac{\partial Z}{\partial t}(\bar{x},\bar{y},\bar{t}) 
&= \frac{1}{2} \, \frac{\partial \Phi}{\partial t}(\bar{x},\bar{t}) \, d(F(\bar{x},\bar{t}),F(\bar{y},\bar{t}))^2 \\ 
&+ H(\bar{x},\bar{t}) + H(\bar{x},\bar{t}) \, \Phi(\bar{x},\bar{t}) \, \big \langle \exp_{F(\bar{x},\bar{t})}^{-1}(F(\bar{y},\bar{t})),\nu(\bar{x},\bar{t}) \big \rangle \\
&- H(\bar{y},\bar{t}) \, \big \langle (D\exp_{F(\bar{x},\bar{t})}^{-1})_{F(\bar{y},\bar{t})}(\nu(\bar{y},\bar{t})),\nu(\bar{x},\bar{t}) + \Phi(\bar{x},\bar{t}) \, \exp_{F(\bar{x},\bar{t})}^{-1}(F(\bar{y},\bar{t})) \big \rangle \\ 
&+ \sum_{i=1}^n \frac{\partial H}{\partial x_i}(\bar{x},\bar{t}) \, \Big \langle \exp_{F(\bar{x},\bar{t})}^{-1}(F(\bar{y},\bar{t})),\frac{\partial F}{\partial x_i}(\bar{x},\bar{t}) \Big \rangle \\ 
&+ H(\bar{x},\bar{t}) \, \Xi_{F(\bar{y},\bar{t}),F(\bar{x},\bar{t})}(\nu(\bar{x},\bar{t}),\nu(\bar{x},\bar{t})). 
\end{align*} 
Consequently, 
\begin{align*}
0 &\geq \frac{\partial Z}{\partial t}(\bar{x},\bar{y},\bar{t}) - \sum_{i=1}^n \Big ( \frac{\partial^2 Z}{\partial x_i^2}(\bar{x},\bar{y},\bar{t}) + 2 \, \frac{\partial^2 Z}{\partial x_i \, \partial y_i}(\bar{x},\bar{y},\bar{t}) + \frac{\partial^2 Z}{\partial y_i^2}(\bar{x},\bar{y},\bar{t}) \Big ) \\ 
&\geq \frac{1}{2} \, \bigg ( \frac{\partial \Phi}{\partial t}(\bar{x},\bar{t}) - \Delta \Phi(\bar{x},\bar{t}) - |A(\bar{x},\bar{t})|^2 \, \Phi(\bar{x},\bar{t}) \\ 
&\hspace{20mm} + \sum_{i=1}^n \frac{2}{\Phi(\bar{x},\bar{t})-\lambda_i(\bar{x},\bar{t})} \, \Big ( \frac{\partial \Phi}{\partial x_i}(\bar{x},\bar{t}) \Big )^2 \bigg ) \, d(F(\bar{x},\bar{t}),F(\bar{y},\bar{t}))^2 \\ 
&- O \big ( d(F(\bar{x},\bar{t}),F(\bar{y},\bar{t})) + H(\bar{x},\bar{t}) \, d(F(\bar{x},\bar{t}),F(\bar{y},\bar{t}))^2 \big ) \\ 
&- O \bigg ( \sum_{i=1}^n \frac{1}{\Phi(\bar{x},\bar{t})-\lambda_i(\bar{x},\bar{t})} \, \Big | \frac{\partial \Phi}{\partial x_i}(\bar{x},\bar{t}) \Big | \, d(F(\bar{x},\bar{t}),F(\bar{y},\bar{t}))^2 \bigg ). 
\end{align*} 
We now multiply both sides by $\frac{2}{d(F(\bar{x},\bar{t}),F(\bar{y},\bar{t}))^2}$. Using the estimate 
\begin{align*} 
\frac{1}{d(F(\bar{x},\bar{t}),F(\bar{y},\bar{t}))} 
&\leq \frac{|\langle \exp_{F(\bar{x},\bar{t})}^{-1}(F(\bar{y},\bar{t})),\nu(\bar{x},\bar{t}) \rangle|}{d(F(\bar{x},\bar{t}),F(\bar{y},\bar{t}))^2} + \sum_{i=1}^n \frac{|\langle \exp_{F(\bar{x},\bar{t})}^{-1}(F(\bar{y},\bar{t})),\frac{\partial F}{\partial x_i}(\bar{x},\bar{t}) \rangle|}{d(F(\bar{x},\bar{t}),F(\bar{y},\bar{t}))^2} \\ 
&\leq \frac{1}{2} \, \Phi(\bar{x},\bar{t}) + \sum_{i=1}^n \frac{1}{2} \, \frac{1}{\Phi(\bar{x},\bar{t}) - \lambda_i(\bar{x},\bar{t})} \, \Big ( \Big | \frac{\partial \Phi}{\partial x_i}(\bar{x},\bar{t}) \Big | + O(1) \Big ), 
\end{align*} 
we obtain  
\begin{align*}
&\frac{\partial \Phi}{\partial t}(\bar{x},\bar{t}) - \Delta \Phi(\bar{x},\bar{t}) - |A(\bar{x},\bar{t})|^2 \, \Phi(\bar{x},\bar{t}) + \sum_{i=1}^n \frac{2}{\Phi(\bar{x},\bar{t})-\lambda_i(\bar{x},\bar{t})} \, \Big ( \frac{\partial \Phi}{\partial x_i}(\bar{x},\bar{t}) \Big )^2 \\ 
&\leq O \bigg ( H(\bar{x},\bar{t}) \, \Phi(\bar{x},\bar{t}) + \sum_{i=1}^n \frac{1}{\Phi(\bar{x},\bar{t}) - \lambda_i(\bar{x},\bar{t})} + \sum_{i=1}^n \frac{1}{\Phi(\bar{x},\bar{t}) - \lambda_i(\bar{x},\bar{t})} \, \Big | \frac{\partial \Phi}{\partial x_i}(\bar{x},\bar{t}) \Big | \bigg ). 
\end{align*} From this, the assertion follows. \\

\begin{corollary} 
\label{evolution.of.mu}
The function $\mu$ satisfies the evolution equation 
\[\frac{\partial \mu}{\partial t} - \Delta \mu - |A|^2 \, \mu + \sum_{i=1}^n \frac{1}{\mu-\lambda_i} \, (D_i \mu)^2 \leq C \, H + C \, \mu + C \, \sum_{i=1}^n \frac{1}{\mu - \lambda_i}\] 
on the set $\{\lambda_n < \mu\} \cap \{\mu \geq 8 \, \text{\rm inj}(N)^{-1}\}$. Here, $\Delta \mu$ is interpreted in the sense of distributions. Moreover, $C$ is a positive constant that depends only on the ambient manifold $N$ and the initial hypersurface $M_0$.
\end{corollary}

\begin{corollary} 
\label{maximum.principle.bound.for.mu}
We have 
\[\sup_{t \in [0,T)} \sup_{M_t} \frac{\mu}{H} \leq C,\] 
where $C$ is a constant that depends only on the ambient manifold $N$, the initial hypersurface $M_0$, and on $T$.
\end{corollary} 

\textbf{Proof.} 
The ratio $\frac{\mu}{H}$ saisfies an evolution equation of the form 
\[\frac{\partial}{\partial t} \Big ( \frac{\mu}{H} \Big ) - \Delta \Big ( \frac{\mu}{H} \Big ) - 2 \, \Big \langle \frac{\nabla H}{H},\nabla \Big ( \frac{\mu}{H} \Big ) \Big \rangle \leq C + C \, \frac{\mu}{H} + C \, \sum_{i=1}^n \frac{1}{H \, (\mu - \lambda_i)}.\] 
It follows from results in \cite{Huisken-Sinestrari2} that 
\[\sup_{t \in [0,T)} \sup_{M_t} \frac{|\lambda_i|+1}{H} \leq K,\] 
where $K$ is a constant that depends only on the ambient manifold $N$, the initial hypersurface $M_0$, and on $T$. Hence, if $\frac{\mu}{H} \geq 2K$, then $\frac{\mu-\lambda_i}{H} \geq K$, and therefore $\frac{1}{H \, (\mu-\lambda_i)} \leq \frac{1}{K \, H^2} \leq K$. Thus, we conclude that 
\[\frac{\partial}{\partial t} \Big ( \frac{\mu}{H} \Big ) - \Delta \Big ( \frac{\mu}{H} \Big ) - 2 \, \Big \langle \frac{\nabla H}{H},\nabla \Big ( \frac{\mu}{H} \Big ) \Big \rangle \leq C + C \, \frac{\mu}{H}\] 
whenever $\frac{\mu}{H} \geq 2K$. Hence, the assertion follows from the maximum principle. \\

Corollary \ref{maximum.principle.bound.for.mu} generalizes the noncollapsing estimate of Andrews \cite{Andrews} to Riemannian manifolds.

\section{An auxiliary inequality}

In this section, we will consider a single hypersurface $M_{\bar{t}}$ for some fixed time $\bar{t}$. We will suppress $\bar{t}$ in the notation, as we will only work with a fixed hypersurface. By the convexity estimate of Huisken and Sinestrari \cite{Huisken-Sinestrari2}, we have a pointwise estimate of the form $\lambda_1 \geq -\varepsilon \, H - K_1(\varepsilon)$, where $\varepsilon$ is an arbitrary positive real number.

\begin{proposition}
\label{aux}
Consider a point $\bar{x} \in M$ such that $\lambda_n(\bar{x}) < \mu(\bar{x})$ and $\mu(\bar{x}) \geq 8 \, \text{\rm inj}(N)^{-1}$. Furthermore, we assume that $U$ is an open neighborhood of $\bar{x}$ and $\Phi: U \to \mathbb{R}$ is a smooth function such that $\Phi(\bar{x}) = \mu(\bar{x})$ and $\Phi(x) \geq \mu(x)$ for all $x \in U$. Then 
\begin{align*} 
0 &\leq \Delta \Phi + \frac{1}{2} \, |A|^2 \, \Phi - \frac{1}{2} \, H \, \Phi^2 + \frac{1}{2} \, n^3 \, (n\varepsilon \, \Phi + K_1(\varepsilon)) \, \Phi^2 \\
&+ \sum_{i=1}^n \frac{1}{\Phi-\lambda_i} \, (|D_i \Phi| + C) \, |D_i H| \\
&+ \big ( H + n^3 \, (n\varepsilon \, \Phi + K_1(\varepsilon)) \big ) \, \sum_{i=1}^n \frac{1}{(\Phi-\lambda_i)^2} \, ((D_i \Phi)^2 + C)^2 \\ 
&+ C \, \Phi + C \, \sum_{i=1}^n \frac{1}{\Phi-\lambda_i} 
\end{align*}
at the point $\bar{x}$. Here, $C$ is a positive constant that depends only on $N$, $M_0$, and $T$.
\end{proposition}

\textbf{Proof.} 
As above, we define
\[Z(x,y) = \frac{1}{2} \, \Phi(x) \, d(F(x),F(y))^2 + \langle \exp_{F(x)}^{-1}(F(y)),\nu(x) \rangle.\] 
By assumption, we have $Z(x,y) \geq 0$ whenever $x \in U$ and $d(F(x),F(y)) \leq \frac{1}{2} \, \text{\rm inj}(N)$. Moreover, there exists a point $\bar{y} \in M$ such that $0 < d(F(\bar{x}),F(\bar{y})) \leq \frac{1}{2} \, \text{\rm inj}(N)$ and $Z(\bar{x},\bar{y}) = 0$. As above, it is easy to see that $\Phi(\bar{x}) \, d(F(\bar{x}),F(\bar{y})) \leq 2$, so $d(F(\bar{x}),F(\bar{y})) \leq \frac{1}{4} \, \text{\rm inj}(N)$. Moreover, we have $H(\bar{x}) \leq C \, \Phi(\bar{x})$ and $H(\bar{y}) \leq C \, \Phi(\bar{x})$ for some constant $C$ that depends only on the ambient manifold $N$.

It follows from results in Section \ref{calc} that 
\begin{align*}
&\sum_{i=1}^n \frac{\partial^2 Z}{\partial x_i^2}(\bar{x},\bar{y}) \\
&\leq \frac{1}{2} \, \bigg ( \Delta \Phi(\bar{x}) + |A(\bar{x})|^2 \, \Phi(\bar{x}) - H(\bar{x}) \, \Phi(\bar{x})^2 \\
&\hspace{20mm} - \sum_{i=1}^n \frac{2}{\Phi(\bar{x})-\lambda_i(\bar{x})} \, \Big ( \frac{\partial \Phi}{\partial x_i}(\bar{x}) \Big )^2 \\ 
&\hspace{20mm} + \sum_{i=1}^n \frac{1}{\Phi(\bar{x})-\lambda_i(\bar{x})} \, \Big ( \frac{\partial \Phi}{\partial x_i}(\bar{x}) + O(1) \Big ) \, \frac{\partial H}{\partial x_i}(\bar{x}) \bigg ) \, d(F(\bar{x}),F(\bar{y}))^2 \\
&+ n \, \Phi(\bar{x}) - H(\bar{x}) \\ 
&+ O \big ( d(F(\bar{x}),F(\bar{y})) \big ) + O \bigg ( \sum_{i=1}^n \frac{1}{\Phi(\bar{x})-\lambda_i(\bar{x})} \, \Big | \frac{\partial \Phi}{\partial x_i}(\bar{x}) \Big | \, d(F(\bar{x}),F(\bar{y}))^2 \bigg ). 
\end{align*}
Moreover, we have
\[\frac{\partial^2 Z}{\partial x_i \, \partial y_i}(\bar{x},\bar{y}) = -(\Phi(\bar{x})-\lambda_i(\bar{x})) + O \big ( d(F(\bar{x},\bar{t}),F(\bar{y},\bar{t})) \big )\]
and
\[\frac{\partial^2 Z}{\partial y_i^2}(\bar{x},\bar{y}) = \Phi(\bar{x}) - h_{ii}(\bar{y}) + O \big ( d(F(\bar{x},\bar{t}),F(\bar{y},\bar{t})) \big ).\]
In particular, we have $h_{ii}(\bar{y}) \leq \Phi(\bar{x}) + O \big ( d(F(\bar{x},\bar{t}),F(\bar{y},\bar{t})) \big )$, hence $H(\bar{y}) \leq n \, \Phi(\bar{x}) + O \big ( d(F(\bar{x},\bar{t}),F(\bar{y},\bar{t})) \big )$. Consequently, the convexity estimate of Huisken and Sinestrari \cite{Huisken-Sinestrari2} implies that $h_{ii}(\bar{y}) \geq -\varepsilon \, H(\bar{y}) - K_1(\varepsilon) \geq -n\varepsilon \, \Phi(\bar{x}) - K_1(\varepsilon) - O \big ( d(F(\bar{x},\bar{t}),F(\bar{y},\bar{t})) \big )$. From this, we deduce that
\[\frac{\partial^2 Z}{\partial y_i^2}(\bar{x},\bar{y}) \leq \Phi(\bar{x}) + n\varepsilon \, \Phi(\bar{x}) + K_1(\varepsilon) + O \big ( d(F(\bar{x},\bar{t}),F(\bar{y},\bar{t})) \big ).\]
Thus, we conclude that
\begin{align*}
&\sum_{i=1}^n \Big ( \frac{\partial^2 Z}{\partial x_i^2}(\bar{x},\bar{y}) + 2 \, \frac{\Phi(\bar{x})-\lambda_i(\bar{x})}{\Phi(\bar{x})} \, \frac{\partial^2 Z}{\partial x_i \, \partial y_i}(\bar{x},\bar{y}) \\
&\hspace{20mm} + \frac{(\Phi(\bar{x})-\lambda_i(\bar{x}))^2}{\Phi(\bar{x})^2} \, \frac{\partial^2 Z}{\partial y_i^2}(\bar{x},\bar{y}) \Big ) \\
&\leq \frac{1}{2} \, \bigg ( \Delta \Phi(\bar{x}) + |A(\bar{x})|^2 \, \Phi(\bar{x}) - H(\bar{x}) \, \Phi(\bar{x})^2 \\
&\hspace{20mm} - \sum_{i=1}^n \frac{2}{\Phi(\bar{x})-\lambda_i(\bar{x})} \, \Big ( \frac{\partial \Phi}{\partial x_i}(\bar{x}) \Big )^2 \\ 
&\hspace{20mm} + \sum_{i=1}^n \frac{1}{\Phi(\bar{x})-\lambda_i(\bar{x})} \, \Big ( \frac{\partial \Phi}{\partial x_i}(\bar{x}) + O(1) \Big ) \, \frac{\partial H}{\partial x_i}(\bar{x}) \bigg ) \, d(F(\bar{x}),F(\bar{y}))^2 \\
&+ n \, \Phi(\bar{x}) - H(\bar{x}) - \sum_{i=1}^n \frac{(\Phi(\bar{x})-\lambda_i(\bar{x}))^2}{\Phi(\bar{x})} \\
&+ \sum_{i=1}^n \frac{(\Phi(\bar{x})-\lambda_i(\bar{x}))^2}{\Phi(\bar{x})^2} \, (n\varepsilon \, \Phi(\bar{x}) + K_1(\varepsilon)) \\ 
&+ O \big ( d(F(\bar{x},\bar{t}),F(\bar{y},\bar{t})) \big ) + O \bigg ( \sum_{i=1}^n \frac{1}{\Phi(\bar{x},\bar{t})-\lambda_i(\bar{x},\bar{t})} \, \Big | \frac{\partial \Phi}{\partial x_i}(\bar{x},\bar{t}) \Big | \, d(F(\bar{x},\bar{t}),F(\bar{y},\bar{t}))^2 \bigg ) \\ 
&\leq \frac{1}{2} \, \bigg ( \Delta \Phi(\bar{x}) + |A(\bar{x})|^2 \, \Phi(\bar{x}) - H(\bar{x}) \, \Phi(\bar{x})^2 \\
&\hspace{20mm} - \sum_{i=1}^n \frac{2}{\Phi(\bar{x})-\lambda_i(\bar{x})} \, \Big ( \frac{\partial \Phi}{\partial x_i}(\bar{x}) \Big )^2 \\ 
&\hspace{20mm} + \sum_{i=1}^n \frac{1}{\Phi(\bar{x})-\lambda_i(\bar{x})} \, \Big ( \frac{\partial \Phi}{\partial x_i}(\bar{x}) + O(1) \Big ) \, \frac{\partial H}{\partial x_i}(\bar{x}) \bigg ) \, d(F(\bar{x}),F(\bar{y}))^2 \\
&+ H(\bar{x}) - \frac{|A(\bar{x})|^2}{\Phi(\bar{x})} + n^3 \, (n\varepsilon \, \Phi(\bar{x}) + K_1(\varepsilon)) \\ 
&+ O \, \big ( d(F(\bar{x},\bar{t}),F(\bar{y},\bar{t})) \big ) + O \bigg ( \sum_{i=1}^n \frac{1}{\Phi(\bar{x},\bar{t})-\lambda_i(\bar{x},\bar{t})} \, \Big | \frac{\partial \Phi}{\partial x_i}(\bar{x},\bar{t}) \Big | \, d(F(\bar{x},\bar{t}),F(\bar{y},\bar{t}))^2 \bigg ). 
\end{align*}
In the last step, we have used the fact that $0 \leq \Phi(\bar{x}) - \lambda_i(\bar{x}) \leq n \, \Phi(\bar{x})$ for $i=1,\hdots,n$ and $\sum_{i=1}^n \frac{(\Phi(\bar{x})-\lambda_i(\bar{x}))^2}{\Phi(\bar{x})^2} \leq n^3$. We now multiply both sides by $\frac{2}{d(F(\bar{x}),F(\bar{y}))^2}$. Using the identity
\begin{align*}
\frac{1}{d(F(\bar{x}),F(\bar{y}))^2} &= \frac{\langle \exp_{F(\bar{x})}^{-1}(F(\bar{y})),\nu(\bar{x}) \rangle^2}{d(F(\bar{x}),F(\bar{y}))^4} + \sum_{i=1}^n \frac{\langle \exp_{F(\bar{x})}^{-1}(F(\bar{y})),\frac{\partial F}{\partial x_i}(\bar{x}) \rangle^2}{d(F(\bar{x}),F(\bar{y}))^4} \\
&= \frac{1}{4} \, \bigg ( \Phi(\bar{x})^2 + \sum_{i=1}^n \frac{1}{(\Phi(\bar{x})-\lambda_i(\bar{x}))^2} \, \Big ( \frac{\partial \Phi}{\partial x_i}(\bar{x}) + O(1) \Big )^2 \bigg ),
\end{align*}
we derive the estimate
\begin{align*}
&\frac{2}{d(F(\bar{x}),F(\bar{y}))^2} \, \sum_{i=1}^n \Big ( \frac{\partial^2 Z}{\partial x_i^2}(\bar{x},\bar{y}) + 2 \, \frac{\Phi(\bar{x})-\lambda_i(\bar{x})}{\Phi(\bar{x})} \, \frac{\partial^2 Z}{\partial x_i \, \partial y_i}(\bar{x},\bar{y}) \\
&\hspace{40mm} + \frac{(\Phi(\bar{x})-\lambda_i(\bar{x}))^2}{\Phi(\bar{x})^2} \, \frac{\partial^2 Z}{\partial y_i^2}(\bar{x},\bar{y}) \Big ) \\
&\leq \Delta \Phi(\bar{x}) + |A(\bar{x})|^2 \, \Phi(\bar{x}) - H(\bar{x}) \, \Phi(\bar{x})^2 \\ 
&- \sum_{i=1}^n \frac{2}{\Phi(\bar{x})-\lambda_i(\bar{x})} \, \Big ( \frac{\partial \Phi}{\partial x_i}(\bar{x}) \Big )^2 + \sum_{i=1}^n \frac{1}{\Phi(\bar{x})-\lambda_i(\bar{x})} \, \Big ( \frac{\partial \Phi}{\partial x_i}(\bar{x}) + O(1) \Big ) \, \frac{\partial H}{\partial x_i}(\bar{x}) \\
&+ \frac{1}{2} \, \Big ( H(\bar{x}) - \frac{|A(\bar{x})|^2}{\Phi(\bar{x})} + n^3 \, (n\varepsilon \, \Phi(\bar{x}) + K_1(\varepsilon)) \Big ) \\
&\hspace{20mm} \cdot \bigg ( \Phi(\bar{x})^2 + \sum_{i=1}^n \frac{1}{(\Phi(\bar{x})-\lambda_i(\bar{x}))^2} \, \Big ( \frac{\partial \Phi}{\partial x_i}(\bar{x}) + O(1) \Big )^2 \bigg ) \\
&+ O \bigg ( \Phi(\bar{x}) + \sum_{i=1}^n \frac{1}{\Phi(\bar{x})-\lambda_i(\bar{x})} + \sum_{i=1}^n \frac{1}{\Phi(\bar{x})-\lambda_i(\bar{x})} \, \Big | \frac{\partial \Phi}{\partial x_i}(\bar{x}) \Big | \bigg ) \\ 
&= \Delta \Phi(\bar{x}) + \frac{1}{2} \, |A(\bar{x})|^2 \, \Phi(\bar{x}) - \frac{1}{2} \, H(\bar{x}) \, \Phi(\bar{x})^2 + \frac{1}{2} \, n^3 \, (n\varepsilon \, \Phi(\bar{x}) + K_1(\varepsilon)) \, \Phi(\bar{x})^2 \\
&- \sum_{i=1}^n \frac{2}{\Phi(\bar{x})-\lambda_i(\bar{x})} \, \Big ( \frac{\partial \Phi}{\partial x_i}(\bar{x}) \Big )^2 + \sum_{i=1}^n \frac{1}{\Phi(\bar{x})-\lambda_i(\bar{x})} \, \Big ( \frac{\partial \Phi}{\partial x_i}(\bar{x}) + O(1) \Big ) \, \frac{\partial H}{\partial x_i}(\bar{x}) \\
&+ \frac{1}{2} \, \Big ( H(\bar{x}) - \frac{|A(\bar{x})|^2}{\Phi(\bar{x})} + n^3 \, (n\varepsilon \, \Phi(\bar{x}) + K_1(\varepsilon)) \Big ) \, \sum_{i=1}^n \frac{1}{(\Phi(\bar{x})-\lambda_i(\bar{x}))^2} \, \Big ( \frac{\partial \Phi}{\partial x_i}(\bar{x}) + O(1) \Big )^2 \\ 
&+ O \bigg ( \Phi(\bar{x}) + \sum_{i=1}^n \frac{1}{\Phi(\bar{x})-\lambda_i(\bar{x})} + \sum_{i=1}^n \frac{1}{\Phi(\bar{x})-\lambda_i(\bar{x})} \, \Big | \frac{\partial \Phi}{\partial x_i}(\bar{x}) \Big | \bigg ). 
\end{align*}
Since the function $Z$ attains a local minimum at the point $(\bar{x},\bar{y})$, we have
\[\sum_{i=1}^n \Big ( \frac{\partial^2 Z}{\partial x_i^2}(\bar{x},\bar{y}) + 2 \, \frac{\Phi(\bar{x})-\lambda_i(\bar{x})}{\Phi(\bar{x})} \, \frac{\partial^2 Z}{\partial x_i \, \partial y_i}(\bar{x},\bar{y}) + \frac{(\Phi(\bar{x})-\lambda_i(\bar{x}))^2}{\Phi(\bar{x})^2} \, \frac{\partial^2 Z}{\partial y_i^2}(\bar{x},\bar{y}) \Big ) \geq 0.\]
Putting these facts together, we obtain 
\begin{align*} 
0 &\leq \Delta \Phi(\bar{x}) + \frac{1}{2} \, |A(\bar{x})|^2 \, \Phi(\bar{x}) - \frac{1}{2} \, H(\bar{x}) \, \Phi(\bar{x})^2 + \frac{1}{2} \, n^3 \, (n\varepsilon \, \Phi(\bar{x}) + K_1(\varepsilon)) \, \Phi(\bar{x})^2 \\
&- \sum_{i=1}^n \frac{2}{\Phi(\bar{x})-\lambda_i(\bar{x})} \, \Big ( \frac{\partial \Phi}{\partial x_i}(\bar{x}) \Big )^2 + \sum_{i=1}^n \frac{1}{\Phi(\bar{x})-\lambda_i(\bar{x})} \, \Big ( \frac{\partial \Phi}{\partial x_i}(\bar{x}) + O(1) \Big ) \, \frac{\partial H}{\partial x_i}(\bar{x}) \\
&+ \frac{1}{2} \, \Big ( H(\bar{x}) - \frac{|A(\bar{x})|^2}{\Phi(\bar{x})} + n^3 \, (n\varepsilon \, \Phi(\bar{x}) + K_1(\varepsilon)) \Big ) \, \sum_{i=1}^n \frac{1}{(\Phi(\bar{x})-\lambda_i(\bar{x}))^2} \, \Big ( \frac{\partial \Phi}{\partial x_i}(\bar{x}) + O(1) \Big )^2 \\ 
&+ O \bigg ( \Phi(\bar{x}) + \sum_{i=1}^n \frac{1}{\Phi(\bar{x})-\lambda_i(\bar{x})} + \sum_{i=1}^n \frac{1}{\Phi(\bar{x})-\lambda_i(\bar{x})} \, \Big | \frac{\partial \Phi}{\partial x_i}(\bar{x}) \Big | \bigg ). 
\end{align*} From this, the assertion follows. \\

\begin{corollary} 
We have 
\begin{align*} 
0 &\leq \Delta \mu + \frac{1}{2} \, |A|^2 \, \mu - \frac{1}{2} \, H \, \mu^2 + \frac{1}{2} \, n^3 \, (n\varepsilon \, \mu + K_1(\varepsilon)) \, \mu^2 \\ 
&+ \sum_{i=1}^n \frac{1}{\mu-\lambda_i} \, (|D_i \mu| + C) \, |D_i H| \\ 
&+ \big ( H + n^3 \, (n\varepsilon \, \mu + K_1(\varepsilon)) \big ) \, \sum_{i=1}^n \frac{1}{(\mu-\lambda_i)^2} \, ((D_i \mu)^2 + C) \\ 
&+ C \, \mu + C \, \sum_{i=1}^n \frac{1}{\mu - \lambda_i} 
\end{align*} 
on the set $\{\lambda_n < \mu\} \cap \{\mu \geq 8 \, \text{\rm inj}(N)^{-1}\}$. Here, $\Delta \mu$ is interpreted in the sense of distributions.
\end{corollary}

\begin{corollary} 
\label{distributional.version.of.aux}
We have 
\begin{align*} 
0 &\leq -\int_{M_t} \langle \nabla \eta,\nabla \mu \rangle + \frac{1}{2} \int_{M_t} \eta \, \big ( |A|^2 \, \mu - H \, \mu^2 + n^3 \, (n\varepsilon \, \mu + K_1(\varepsilon)) \, \mu^2 \big ) \\ 
&+ \int_{M_t} \eta \, \sum_{i=1}^n \frac{1}{\mu-\lambda_i} \, (|D_i \mu| + C) \, |D_i H| \\ 
&+ \int_{M_t} \eta \, \big ( H + n^3 \, (n\varepsilon \, \mu + K_1(\varepsilon)) \big ) \, \sum_{i=1}^n \frac{1}{(\mu-\lambda_i)^2} \, ((D_i \mu)^2 + C) \\ 
&+ C \int_{M_t} \eta \, \mu + C \int_{M_t} \eta \, \sum_{i=1}^n \frac{1}{\mu - \lambda_i} 
\end{align*} 
for every nonnegative test function $\eta$ which is supported in the set $\{\lambda_n < \mu\} \cap \{\mu \geq 8 \, \text{\rm inj}(N)^{-1}\}$.
\end{corollary}

\section{Proof of Theorem \ref{inscribed.radius}}

Let us fix positive real numbers $\delta>0$ and $T>0$. By the convexity estimate of Huisken and Sinestrari \cite{Huisken-Sinestrari2}, we can find a constant $K_0 \geq 8 \, \text{\rm inj}(N)^{-1} \, \big ( \inf_{t \in [0,T)} \inf_{M_t} \min\{H,1\} \big )^{-1}$ such that
\[(n-1) \, \lambda_1 \geq -\frac{\delta}{2} \, H - K_0 \, \min \{H,1\}\] 
for $t \in [0,T)$. Here, $K_0$ is a constant that depends only on $N$, $M_0$, $\delta$, and $T$.
 
For each $\sigma \in (0,\frac{1}{2})$, we define
\[f_\sigma = H^{\sigma-1} \, (\mu - (1+\delta) \, H) - K_0\]
and
\[f_{\sigma,+} = \max \{f_\sigma,0\}.\]
On the set $\{f_\sigma \geq 0\}$, we have
\[\mu \geq (1+\delta) \, H + K_0 \, H^{1-\sigma} \geq (1+\delta) \, H + K_0 \, \min \{H,1\},\]
hence
\[\mu - \lambda_n \geq \sum_{i=1}^{n-1} \lambda_i + \delta \, H + K_0 \, \min \{H,1\} \geq \frac{\delta}{2} \, H.\]
In particular, we have $\{f_\sigma \geq 0\} \subset \{\lambda_n < \mu\} \cap \{\mu \geq 8 \, \text{\rm inj}(N)^{-1}\}$. By Corollary \ref{maximum.principle.bound.for.mu}, we can find a constant $\Lambda \geq 1$, depending only on $N$, $M_0$, and $T$, such that $\mu \leq \Lambda \, H$ and $|A|^2 \leq \Lambda \, H^2$ for $t \in [0,T)$.

\begin{proposition}
\label{Lp.bound}
Given any $\delta>0$, we can find a positive constant $c_0$, depending only on $\delta$ and the initial hypersurface $M_0$, with the following property: if $p \geq \frac{1}{c_0}$ and $\sigma \leq c_0 \, p^{-\frac{1}{2}}$, then we have
\[\frac{d}{dt} \bigg ( \int_{M_t} f_{\sigma,+}^p \bigg ) \leq C \, \sigma \, p \int_{M_t} f_{\sigma,+}^p + \sigma \, p \, K_0^p \int_{M_t} |A|^2 + C \, p^p \int_{M_t} H^{2-(2-\sigma) \, p}\]
for almost all $t \in [0,T)$. Here, $C$ is a positive constant that depends only on $N$, $M_0$, $\delta$, and $T$, but not on $\sigma$ and $p$.
\end{proposition}

\textbf{Proof.}
By Corollary \ref{evolution.of.mu}, we have 
\[\frac{\partial \mu}{\partial t} - \Delta \mu - |A|^2 \, \mu + \sum_{i=1}^n \frac{1}{\mu-\lambda_i} \, (D_i \mu)^2 \leq C \, H\]
on the set $\{f_\sigma \geq 0\}$, where $\Delta \mu$ is interpreted in the sense of distributions, and $C$ is a positive constant that depends only on $N$, $M_0$, $\delta$, and $T$. A straightforward calculation gives
\begin{align*}
&\frac{\partial}{\partial t} f_\sigma - \Delta f_\sigma - 2 \, (1-\sigma) \, \Big \langle \frac{\nabla H}{H},\nabla f_\sigma \Big \rangle \\
&+ \sum_{i=1}^n \frac{H^{\sigma-1}}{\mu-\lambda_i} \, (D_i \mu)^2 - \sigma \, |A|^2 \, (f_\sigma + K_0) \\
&\leq -\sigma \, (1-\sigma) \, H^{\sigma-3} \, (\mu - (1+\delta) \, H) \, |\nabla H|^2 + C \, H^\sigma \\ 
&\leq C \, H^\sigma
\end{align*}
on the set $\{f_\sigma \geq 0\}$, where $\Delta f_\sigma$ is again interpreted in the sense of distributions. This implies
\begin{align*}
&\frac{d}{dt} \bigg ( \int_{M_t} f_{\sigma,+}^p \bigg ) \\
&\leq -p(p-1) \int_{M_t} f_{\sigma,+}^{p-2} \, |\nabla f_\sigma|^2 + 2 \, (1-\sigma) \, p \int_{M_t} f_{\sigma,+}^{p-1} \, \Big \langle \frac{\nabla H}{H},\nabla f_\sigma \Big \rangle \\
&- p \int_{M_t} \sum_{i=1}^n \frac{f_{\sigma,+}^{p-1} \, H^{\sigma-1}}{\mu-\lambda_i} \, (D_i \mu)^2 + \sigma \, p \int_{M_t} |A|^2 \, f_{\sigma,+}^{p-1} \, (f_\sigma+K_0) \\ 
&+ \int_{M_t} (C \, p \, H^\sigma \, f_{\sigma,+}^{p-1} - H^2 \, f_{\sigma,+}^p). 
\end{align*}
The integral of $|A|^2 \, f_{\sigma,+}^{p-1} \, (f_\sigma+K_0)$ has an unfavorable sign. To estimate this term, we put $\varepsilon = \frac{\delta}{4n^4\Lambda^2}$. Applying Corollary \ref{distributional.version.of.aux} to the test function $\eta = \frac{f_{\sigma,+}^p}{H}$ gives
\begin{align*}
&\frac{1}{2} \int_{M_t} \big ( H \, \mu^2 - |A|^2 \, \mu - n^3 \, (n\varepsilon \, \mu + K_1(\varepsilon)) \, \mu^2 \big ) \, \frac{f_{\sigma,+}^p}{H} \\
&\leq -\int_M \Big \langle \nabla \Big ( \frac{f_{\sigma,+}^p}{H} \Big ),\nabla \mu \Big \rangle + \int_M \frac{f_{\sigma,+}^p}{H} \, \sum_{i=1}^n \frac{1}{\mu-\lambda_i} \, (|D_i \mu|+C) \, |D_i H| \\
&+ \int_M \frac{f_{\sigma,+}^p}{H} \, \big ( H + n^3 \, (n\varepsilon \, \mu + K_1(\varepsilon)) \big ) \, \sum_{i=1}^n \frac{1}{(\mu-\lambda_i)^2} \, ((D_i \mu)^2+C) \\ 
&+ C \int_{M_t} \frac{f_{\sigma,+}^p}{H} \, \mu + C \int_{M_t} \frac{f_{\sigma,+}^p}{H} \, \sum_{i=1}^n \frac{1}{\mu-\lambda_i} \\ 
&\leq p \int_{M_t} \frac{f_{\sigma,+}^{p-1}}{H} \, |\nabla \mu| \, |\nabla f_\sigma| + C \int_{M_t} \frac{f_{\sigma,+}^p}{H^2} \, (|\nabla \mu|+1) \, |\nabla H| \\ 
&+ C \int_{M_t} \frac{f_{\sigma,+}^p}{H^2} \, |\nabla \mu|^2 + C \int_{M_t} f_{\sigma,+}^p. 
\end{align*}
Here, $C$ is a positive constant which depends on $N$, $M_0$, $\delta$, and $T$, but not on $\sigma$ and $p$. On the set $\{f_\sigma \geq 0\}$, we have $\mu \geq (1+\delta) \, H$. Moreover, the convexity estimate of Huisken and Sinestrari implies that $|A|^2 \leq (1+\varepsilon) \, H^2 + K_2(\varepsilon)$. Consequently, we have
\begin{align*}
&H \, \mu^2 - |A|^2 \, \mu - n^3 \, (n\varepsilon \, \mu + K_1(\varepsilon)) \, \mu^2 \\
&\geq (1+\delta) \, H^2 \, \mu - |A|^2 \, \mu - n^3 \, (n\varepsilon \, \mu + K_1(\varepsilon)) \, \Lambda^2 \, H^2 \\
&\geq (\delta-\varepsilon) \, H^2 \, \mu - n^3 \, (n\varepsilon \, \mu + K_1(\varepsilon)) \, \Lambda^2 \, H^2 - K_2(\varepsilon) \, \mu \\
&\geq \frac{\delta}{2} \, H^2 \, \mu - C \, H
\end{align*}
on the set $\{f_\sigma \geq 0\}$. Therefore, we obtain
\begin{align*}
\int_{M_t} |A|^2 \, f_{\sigma,+}^p
&\leq C \, p \int_{M_t} \frac{f_{\sigma,+}^{p-1}}{H} \, |\nabla \mu| \, |\nabla f_\sigma| + C \int_{M_t} \frac{f_{\sigma,+}^p}{H^2} \, (|\nabla \mu|+1) \, |\nabla H| \\
&+ C \int_{M_t} \frac{f_{\sigma,+}^p}{H^2} \, |\nabla \mu|^2 + C \int_{M_t} f_{\sigma,+}^p,
\end{align*}
where $C$ is a positive constant that depends only on $N$, $M_0$, $\delta$, and $T$. Using the pointwise inequality
\[f_{\sigma,+}^{p-1} \, (f_\sigma+K_0) \leq 2 \, f_{\sigma,+}^p + K_0^p,\] 
we obtain
\begin{align*}
&\int_{M_t} |A|^2 \, f_{\sigma,+}^{p-1} \, (f_\sigma+K_0) \\
&\leq C \, p \int_{M_t} \frac{f_{\sigma,+}^{p-1}}{H} \, |\nabla \mu| \, |\nabla f_\sigma| + C \int_{M_t} \frac{f_{\sigma,+}^p}{H^2} \, (|\nabla \mu|+1) \, |\nabla H| \\
&+ C \int_{M_t} \frac{f_{\sigma,+}^p}{H^2} \, |\nabla \mu|^2 + C \int_{M_t} f_{\sigma,+}^p + K_0^p \int_{M_t} |A|^2,
\end{align*}
where $C$ is a positive constant that depends only on $N$, $M_0$, $\delta$, and $T$. Putting these facts together, we conclude that
\begin{align*}
&\frac{d}{dt} \bigg ( \int_{M_t} f_{\sigma,+}^p \bigg ) \\
&\leq -p(p-1) \int_{M_t} f_{\sigma,+}^{p-2} \, |\nabla f_\sigma|^2 + 2 \, (1-\sigma) \, p \int_{M_t} f_{\sigma,+}^{p-1} \, \Big \langle \frac{\nabla H}{H},\nabla f_\sigma \Big \rangle \\
&- p \int_{M_t} \sum_{i=1}^n \frac{f_{\sigma,+}^{p-1} \, H^{\sigma-1}}{\mu-\lambda_i} \, (D_i \mu)^2 + C \, \sigma \, p^2 \int_{M_t} \frac{f_{\sigma,+}^{p-1}}{H} \, |\nabla \mu| \, |\nabla f_\sigma| \\
&+ C \, \sigma \, p \int_{M_t} \frac{f_{\sigma,+}^p}{H^2} \, (|\nabla \mu|+1) \, |\nabla H| + C \, \sigma \, p \int_{M_t} \frac{f_{\sigma,+}^p}{H^2} \, |\nabla \mu|^2 \\
&+ C \, \sigma \, p \int_{M_t} f_{\sigma,+}^p + \sigma \, p \, K_0^p \int_{M_t} |A|^2 + \int_{M_t} (C \, p \, H^\sigma \, f_{\sigma,+}^{p-1} - H^2 \, f_{\sigma,+}^p), 
\end{align*}
where $C$ is a positive constant that depends only on $N$, $M_0$, $\delta$, and $T$. Using the identity
\[\frac{\nabla H}{H} = \frac{\nabla \mu - H^{1-\sigma} \, \nabla f_\sigma}{(1-\sigma) \, \mu + \sigma \, (1+\delta) \, H},\]
we obtain
\[\Big \langle \frac{\nabla H}{H},\nabla f_\sigma \Big \rangle \leq \frac{\langle \nabla \mu,\nabla f_\sigma \rangle}{(1-\sigma) \, \mu + \sigma \, (1+\delta) \, H} \leq C \, H^{-1} \, |\nabla \mu| \, |\nabla f_\sigma|\]
and
\[\frac{|\nabla H|}{H} \leq C \, \frac{|\nabla \mu|}{H} + C \, \frac{|\nabla f_\sigma|}{f_{\sigma,+}}.\]
This implies
\begin{align*}
&\frac{d}{dt} \bigg ( \int_{M_t} f_{\sigma,+}^p \bigg ) \\
&\leq -p(p-1) \int_{M_t} f_{\sigma,+}^{p-2} \, |\nabla f_\sigma|^2 - p \int_{M_t} \sum_{i=1}^n \frac{f_{\sigma,+}^{p-1} \, H^{\sigma-1}}{\mu-\lambda_i} \, (D_i \mu)^2 \\
&+ C \, (p + \sigma \, p^2) \int_{M_t} \frac{f_{\sigma,+}^{p-1}}{H} \, |\nabla \mu| \, |\nabla f_\sigma| + C \, \sigma \, p \int_{M_t} \frac{f_{\sigma,+}^p}{H^2} \, |\nabla \mu|^2 \\ 
&+ C \, \sigma \, p \int_{M_t} \frac{f_{\sigma,+}^p}{H^2} \, |\nabla \mu| + C \, \sigma \, p \int_{M_t} \frac{f_{\sigma,+}^{p-1}}{H} \, |\nabla f_\sigma| \\ 
&+ C \, \sigma \, p \int_{M_t} f_{\sigma,+}^p + \sigma \, p \, K_0^p \int_{M_t} |A|^2 + \int_{M_t} (C \, p \, H^\sigma \, f_{\sigma,+}^{p-1} - H^2 \, f_{\sigma,+}^p), 
\end{align*}
where $C$ is a positive constant that depends only on $N$, $M_0$, $\delta$, and $T$. 

Therefore, we can find a positive constant $c_0$, depending only on $\delta$, $M_0$, $N$, and $T$, with the following property: if $p \geq \frac{1}{c_0}$ and $\sigma \leq c_0 \, p^{-\frac{1}{2}}$, then we have
\[\frac{d}{dt} \bigg ( \int_{M_t} f_{\sigma,+}^p \bigg ) \leq C \, \sigma \, p \int_{M_t} f_{\sigma,+}^p + \sigma \, p \, K_0^p \int_{M_t} |A|^2 + \int_{M_t} (C \, p \, H^\sigma \, f_{\sigma,+}^{p-1} - H^2 \, f_{\sigma,+}^p)\]
for almost all $t \in [0,T)$. Since we have a pointwise lower bound for the function $H$ for $t \in [0,T)$, we obtain a pointwise upper bound for the function $C \, p \, H^\sigma \, f_{\sigma,+}^{p-1} - H^2 \, f_{\sigma,+}^p$ for all $t \in [0,T)$. This gives 
\[\frac{d}{dt} \bigg ( \int_{M_t} f_{\sigma,+}^p \bigg ) \leq C \, \sigma \, p \int_{M_t} f_{\sigma,+}^p + \sigma \, p \, K_0^p \int_{M_t} |A|^2 + C \, p^p \int_{M_t} H^{2-(2-\sigma) \, p}\]
for almost all $t \in [0,T)$, where again $C$ depends only on $N$, $M_0$, $\delta$, and $T$. This completes the proof of Proposition \ref{Lp.bound}. \\

As usual, we can now use the Michael-Simon Sobolev inequality (cf. \cite{Michael-Simon}) and Stampacchia iteration to show that 
\[\mu \leq (1+\delta) \, H + C(N,M_0,\delta,T)\] 
for all $t \in [0,T)$ and all points on $M_t$. This is the desired noncollapsing estimate.

\section{Proof of Theorem \ref{outer.radius}}

Finally, we give the proof of Theorem \ref{outer.radius}. 

\begin{proposition}
\label{evolution.of.rho}
The function $\rho$ satisfies 
\[\frac{\partial \rho}{\partial t} - \Delta \rho - |A|^2 \, \rho + \sum_{i=1}^n \frac{1}{\rho+\lambda_i} \, (D_i \rho)^2 \leq C \, H + C \, \rho + C \, \sum_{i=1}^n \frac{1}{\rho + \lambda_i}\]
on the set $\{\rho+\lambda_1 > 0\} \cap \{\rho \geq 8 \, \text{\rm inj}(N)^{-1}\}$. Here, $\Delta \rho$ is interpreted in the sense of distributions. Moreover, $C$ is a positive constant that depends only on the ambient manifold $N$ and the initial hypersurface $M_0$.
\end{proposition}

The proof of Proposition \ref{evolution.of.rho} is analogous to the proof of Corollary \ref{evolution.of.mu} above. In fact, it suffices to reverse the orientation of $M_t$ everywhere in the argument. \\

\begin{corollary}
\label{maximum.principle.bound.for.rho}
We have 
\[\sup_{t \in [0,T)} \sup_{M_t} \frac{\rho}{H} \leq C,\] 
where $C$ is a constant that depends only on the ambient manifold $N$, the initial hypersurface $M_0$, and on $T$.
\end{corollary} 

\textbf{Proof.} 
The ratio $\frac{\rho}{H}$ saisfies an evolution equation of the form 
\[\frac{\partial}{\partial t} \Big ( \frac{\rho}{H} \Big ) - \Delta \Big ( \frac{\rho}{H} \Big ) - 2 \, \Big \langle \frac{\nabla H}{H},\nabla \Big ( \frac{\rho}{H} \Big ) \Big \rangle \leq C + C \, \frac{\rho}{H} + C \, \sum_{i=1}^n \frac{1}{H \, (\rho + \lambda_i)}.\] 
It follows from results in \cite{Huisken-Sinestrari2} that 
\[\sup_{t \in [0,T)} \sup_{M_t} \frac{|\lambda_i|+1}{H} \leq K,\] 
where $K$ is a constant that depends only on the ambient manifold $N$, the initial hypersurface $M_0$, and on $T$. Hence, if $\frac{\rho}{H} \geq 2K$, then $\frac{\rho+\lambda_i}{H} \geq K$, and therefore $\frac{1}{H \, (\rho+\lambda_i)} \leq \frac{1}{K \, H^2} \leq K$. Thus, we conclude that 
\[\frac{\partial}{\partial t} \Big ( \frac{\rho}{H} \Big ) - \Delta \Big ( \frac{\rho}{H} \Big ) - 2 \, \Big \langle \frac{\nabla H}{H},\nabla \Big ( \frac{\rho}{H} \Big ) \Big \rangle \leq C + C \, \frac{\rho}{H}\] 
whenever $\frac{\rho}{H} \geq 2K$. Hence, the assertion follows from the maximum principle. \\

We next state an auxiliary result:

\begin{proposition}
\label{aux.2}
Consider a point $(\bar{x},\bar{t}) \in M \times [0,T)$ such that $\rho(\bar{x},\bar{t}) + \lambda_1(\bar{x},\bar{t}) > 0$ and $\rho(\bar{x},\bar{t}) \geq 8 \, \text{\rm inj}(N)^{-1}$. We further assume that $U \subset M \times [0,T)$ is an open neighborhood of $\bar{x}$ and $\Phi: U \times (\bar{t}-\alpha,\bar{t}] \to \mathbb{R}$ is a smooth function such that $\Phi(\bar{x},\bar{t}) = \rho(\bar{x},\bar{t})$ and $\Phi(x,t) \geq \rho(x,t)$ for all points $(x,t) \in U \times (\bar{t}-\alpha,\bar{t}]$. Then
\begin{align*} 
&\frac{\partial \Phi}{\partial t} + \frac{1}{2} \, H \, \Phi^2 - \sum_{i=1}^n \frac{1}{\Phi+\lambda_i} \, (|D_i \Phi| + L) \, (|D_i H|+L) - \sum_{i=1}^n \frac{H}{(\Phi+\lambda_i)^2} \, ((D_i \Phi)^2 + L) \\ 
&\leq L \,  H 
\end{align*}
at the point $(\bar{x},\bar{t})$. Here, $L$ is a positive constant that depends only on the ambient manifold $N$, the initial hypersurface $M_0$, and on $T$.
\end{proposition}

\textbf{Proof.}
We define
\begin{align*} 
W(x,y,t) 
&= \Phi(x,t) \, \psi_{F(y,t)}(F(x,t)) + \big \langle \nabla \psi_{F(y,t)} \big |_{F(x,t)},\nu(x,t) \big \rangle \\ 
&= \frac{1}{2} \, \Phi(x,t) \, d(F(x,t),F(y,t))^2 - \langle \exp_{F(x,t)}^{-1}(F(y,t)),\nu(x,t) \rangle. 
\end{align*}
By assumption, we have $W(x,y,t) \geq 0$ whenever $x \in U$, $t \in (\bar{t}-\alpha,\bar{t}]$, and $d(F(x,t),F(y,t)) \leq \frac{1}{2} \, \text{\rm inj}(N)$. Moreover, we can find a point $\bar{y}$ such that $0 < d(F(\bar{x},\bar{t}),F(\bar{y},\bar{t})) \leq \frac{1}{2} \, \text{\rm inj}(N)$ and $W(\bar{x},\bar{y},\bar{t}) = 0$. From this, we deduce that $\Phi(\bar{x},\bar{t}) \, d(F(\bar{x},\bar{t}),F(\bar{y},\bar{t})) \leq 2$, hence $d(F(\bar{x},\bar{t}),F(\bar{y},\bar{t})) \leq \frac{1}{4} \, \text{\rm inj}(N)$. As in Section \ref{calc}, we compute 
\begin{align*}
0 = \frac{\partial W}{\partial x_i}(\bar{x},\bar{y},\bar{t}) 
&= \frac{1}{2} \, \frac{\partial \Phi}{\partial x_i}(\bar{x},\bar{t}) \, d(F(\bar{y},\bar{t}),F(\bar{x},\bar{t}))^2 \\ 
&- \Phi(\bar{x},\bar{t}) \, \Big \langle \exp_{F(\bar{x},\bar{t})}^{-1}( F(\bar{y},\bar{t})),\frac{\partial F}{\partial x_i}(\bar{x},\bar{t}) \Big \rangle \\ 
&- h_i^j(\bar{x},\bar{t}) \, \Big \langle \exp_{F(\bar{x},\bar{t})}^{-1}( F(\bar{y},\bar{t})),\frac{\partial F}{\partial x_j}(\bar{x},\bar{t}) \Big \rangle \\ 
&+ \Xi_{F(\bar{y},\bar{t}),F(\bar{x},\bar{t})} \Big ( \frac{\partial F}{\partial x_i}(\bar{x},\bar{t}),\nu(\bar{x},\bar{t}) \Big ). 
\end{align*} 
Let us pick geodesic normal coordinates around $\bar{x}$ such that $h_{ij}(\bar{x},\bar{t})$ is a diagonal matrix. The relation $\frac{\partial W}{\partial x_i}(\bar{x},\bar{y},\bar{t}) = 0$ implies
\begin{align*} 
&\Big \langle \exp_{F(\bar{x},\bar{t})}^{-1}( F(\bar{y},\bar{t})),\frac{\partial F}{\partial x_i}(\bar{x},\bar{t}) \Big \rangle \\ 
&= \frac{1}{2} \, \frac{1}{\Phi(\bar{x},\bar{t}) + \lambda_i(\bar{x},\bar{t})} \, \Big ( \frac{\partial \Phi}{\partial x_i}(\bar{x},\bar{t}) + O(1) \Big ) \, d(F(\bar{x},\bar{t}),F(\bar{x},\bar{t}))^2. 
\end{align*} 
In the next step, we use the identity  
\begin{align*} 
\frac{\partial W}{\partial t}(\bar{x},\bar{y},\bar{t}) 
&= \frac{1}{2} \, \frac{\partial \Phi}{\partial t}(\bar{x},\bar{t}) \, d(F(\bar{x},\bar{t}),F(\bar{y},\bar{t}))^2 \\ 
&- H(\bar{x},\bar{t}) + H(\bar{x},\bar{t}) \, \Phi(\bar{x},\bar{t}) \, \big \langle \exp_{F(\bar{x},\bar{t})}^{-1}(F(\bar{y},\bar{t})),\nu(\bar{x},\bar{t}) \big \rangle \\
&+ H(\bar{y},\bar{t}) \, \big \langle (D\exp_{F(\bar{x},\bar{t})}^{-1})_{F(\bar{y},\bar{t})}(\nu(\bar{y},\bar{t})),\nu(\bar{x},\bar{t}) - \Phi(\bar{x},\bar{t}) \, \exp_{F(\bar{x},\bar{t})}^{-1}(F(\bar{y},\bar{t})) \big \rangle \\ 
&- \sum_{i=1}^n \frac{\partial H}{\partial x_i}(\bar{x},\bar{t}) \, \Big \langle \exp_{F(\bar{x},\bar{t})}^{-1}(F(\bar{y},\bar{t})),\frac{\partial F}{\partial x_i}(\bar{x},\bar{t}) \Big \rangle \\ 
&- H(\bar{x},\bar{t}) \, \Xi_{F(\bar{y},\bar{t}),F(\bar{x},\bar{t})}(\nu(\bar{x},\bar{t}),\nu(\bar{x},\bar{t})). 
\end{align*} 
The terms $H(\bar{y},\bar{t})$ and $\big \langle (D\exp_{F(\bar{x},\bar{t})}^{-1})_{F(\bar{y},\bar{t})}(\nu(\bar{y},\bar{t})),\nu(\bar{x},\bar{t}) - \Phi(\bar{x},\bar{t}) \, \exp_{F(\bar{x},\bar{t})}^{-1}(F(\bar{y},\bar{t})) \big \rangle$ are nonnegative. This gives  
\begin{align*} 
&\frac{\partial W}{\partial t}(\bar{x},\bar{y},\bar{t}) \\ 
&\geq \frac{1}{2} \, \bigg ( \frac{\partial \Phi}{\partial t}(\bar{x},\bar{t}) + H(\bar{x},\bar{t}) \, \Phi(\bar{x},\bar{t})^2 \\
&\hspace{20mm} - \sum_{i=1}^n \frac{1}{\Phi(\bar{x},\bar{t})+\lambda_i(\bar{x},\bar{t})} \, \Big ( \frac{\partial \Phi}{\partial x_i}(\bar{x},\bar{t}) + O(1) \Big ) \, \frac{\partial H}{\partial x_i}(\bar{x},\bar{t}) \bigg ) \, |F(\bar{x},\bar{t}) - F(\bar{y},\bar{t})|^2 \\
&- H(\bar{x},\bar{t}) + O \big ( H(\bar{x},\bar{t}) \, d(F(\bar{x},\bar{t}),F(\bar{y},\bar{t}))^2 \big ). 
\end{align*}
We now multiply both sides by $\frac{2}{|F(\bar{x},\bar{t})-F(\bar{y},\bar{t})|^2}$. Using the relation 
\begin{align*}
\frac{1}{d(F(\bar{x},\bar{t}),F(\bar{y},\bar{t}))^2} &= \frac{\langle \exp_{F(\bar{x},\bar{t})}^{-1}(F(\bar{y},\bar{t})),\nu(\bar{x},\bar{t}) \rangle^2}{d(F(\bar{x},\bar{t}),F(\bar{y},\bar{t}))^4} + \sum_{i=1}^n \frac{\langle \exp_{F(\bar{x},\bar{t})}^{-1}(F(\bar{y},\bar{t})),\frac{\partial F}{\partial x_i}(\bar{x},\bar{t}) \rangle^2}{d(F(\bar{x},\bar{t}),F(\bar{y},\bar{t}))^4} \\
&= \frac{1}{4} \, \bigg ( \Phi(\bar{x},\bar{t})^2 + \sum_{i=1}^n \frac{1}{(\Phi(\bar{x},\bar{t})+\lambda_i(\bar{x}))^2} \, \Big ( \frac{\partial \Phi}{\partial x_i}(\bar{x},\bar{t}) + O(1) \Big )^2 \bigg ),
\end{align*}
we deduce that
\begin{align*}
&\frac{2}{d(F(\bar{x},\bar{t}),F(\bar{y},\bar{t}))^2} \, \frac{\partial W}{\partial t}(\bar{x},\bar{y},\bar{t}) \\
&\geq \frac{\partial \Phi}{\partial t}(\bar{x},\bar{t}) + H(\bar{x},\bar{t}) \, \Phi(\bar{x},\bar{t})^2 - \sum_{i=1}^n \frac{1}{\Phi(\bar{x},\bar{t})+\lambda_i(\bar{x},\bar{t})} \, \Big ( \frac{\partial \Phi}{\partial x_i}(\bar{x},\bar{t}) + O(1) \Big ) \, \frac{\partial H}{\partial x_i}(\bar{x},\bar{t}) \\
&- \frac{1}{2} \, H(\bar{x},\bar{t}) \, \bigg ( \Phi(\bar{x},\bar{t})^2 + \sum_{i=1}^n \frac{1}{(\Phi(\bar{x},\bar{t})+\lambda_i(\bar{x},\bar{t}))^2} \, \Big ( \frac{\partial \Phi}{\partial x_i}(\bar{x},\bar{t}) + O(1) \Big )^2 \bigg ) - O(H(\bar{x},\bar{t})) \\ 
&= \frac{\partial \Phi}{\partial t}(\bar{x},\bar{t}) + \frac{1}{2} \, H(\bar{x},\bar{t}) \, \Phi(\bar{x},\bar{t})^2 - \sum_{i=1}^n \frac{1}{\Phi(\bar{x},\bar{t})+\lambda_i(\bar{x},\bar{t})} \, \Big ( \frac{\partial \Phi}{\partial x_i}(\bar{x},\bar{t}) + O(1) \Big ) \, \frac{\partial H}{\partial x_i}(\bar{x},\bar{t}) \\
&- \frac{1}{2} \, H(\bar{x},\bar{t}) \, \sum_{i=1}^n \frac{1}{(\Phi(\bar{x},\bar{t})+\lambda_i(\bar{x},\bar{t}))^2} \, \Big ( \frac{\partial \Phi}{\partial x_i}(\bar{x},\bar{t}) + O(1) \Big )^2 - O(H(\bar{x},\bar{t})). 
\end{align*}
Since $\frac{\partial W}{\partial t}(\bar{x},\bar{y},\bar{t}) \leq 0$, the assertion follows. \\

\begin{corollary} 
\label{consequence.of.aux.2}
We have
\begin{align*} 
&\frac{\partial \rho}{\partial t} - \sum_{i=1}^n \frac{1}{\rho+\lambda_i} \, (|D_i \rho| + L) \, (|D_i H| + L) - \sum_{i=1}^n \frac{H}{(\rho+\lambda_i)^2} \, ((D_i \rho)^2 + L) \\ 
&\leq L \,  H 
\end{align*}
almost everywhere on the set $\{\rho+\lambda_1 > 0\} \cap \{\rho \geq 8 \, \text{\rm inj}(N)^{-1}\}$.
\end{corollary}

Let $\delta>0$ be given. The convexity estimate of Huisken and Sinestrari \cite{Huisken-Sinestrari2} implies that we can find a constant $K_0 \geq 8 \, \text{\rm inj}(N)^{-1} \, \big ( \inf_{t \in [0,T)} \inf_{M_t} \min\{H,1\} \big )^{-1}$, depending only on $N$, $M_0$, $\delta$, and $T$, such that
\[\lambda_1 \geq -\frac{\delta}{2} \, H - K_0 \, \min \{H,1\}.\]
For each $\sigma \in (0,\frac{1}{2})$, we put
\[g_\sigma = H^{\sigma-1} \, (\rho - \delta \, H) - K_0\]
and
\[g_{\sigma,+} = \max \{g_\sigma,0\}.\]
On the set $\{g_\sigma \geq 0\}$, we have
\[\rho \geq \delta \, H + K_0 \, H^{1-\sigma} \geq \delta \, H + K_0 \, \min \{H,1\},\]
hence
\[\rho + \lambda_1 \geq \frac{\delta}{2} \, H.\]
In particular, we have $\{g_\sigma \geq 0\} \subset \{\rho+\lambda_1 > 0\} \cap \{\rho \geq 8 \, \text{\rm inj}(N)^{-1}\}$. Furthermore, by Corollary \ref{maximum.principle.bound.for.rho}, there exists a constant $\Lambda \geq 1$, depending only on $N$, $M_0$, and $T$, such that $\rho \leq \Lambda \, H$ and $|A|^2 \leq \Lambda \, H^2$ for $t \in [0,T)$.

\begin{proposition}
\label{Lp.bound.2}
Given any $\delta>0$, there exists a positive constant $c_0$, depending only on $N$, $M_0$, $\delta$, and $T$, with the following property: if $p \geq \frac{1}{c_0}$ and $\sigma \leq c_0 \, p^{-\frac{1}{2}}$, then we have
\[\frac{d}{dt} \bigg ( \int_{M_t} g_{\sigma,+}^p \bigg ) \leq C \, \sigma \, p \int_{M_t} g_{\sigma,+}^p + \sigma \, p \, K_0^p \int_{M_t} |A|^2 + C \, p^p \int_{M_t} H^{2-(2-\sigma)p}\]
for almost all $t$.
\end{proposition}

\textbf{Proof.}
For abbreviation, we define a function $\omega$ by
\begin{align*}
\omega
&= \Delta \rho - \sum_{i=1}^n \frac{1}{\rho+\lambda_i} \, (D_i \rho)^2 \\
&- \sum_i \frac{1}{\rho+\lambda_i} \, (|D_i \rho|+L) \, (|D_i H|+L) - \sum_{i=1}^n \frac{H}{(\rho+\lambda_i)^2} \, ((D_i \rho)^2+L), 
\end{align*}
where $L$ is the constant in Corollary \ref{consequence.of.aux.2}. Combining Proposition \ref{evolution.of.rho} and Corollary \ref{consequence.of.aux.2}, we obtain 
\[\frac{\partial \rho}{\partial t} - \Delta \rho - |A|^2 \, \rho + \sum_{i=1}^n \frac{1}{\rho+\lambda_i} \, (D_i \rho)^2 \leq -\max \{\omega + |A|^2 \, \rho,0\} + C \, H\]
on the set $\{g_\sigma \geq 0\}$. From this, we deduce that
\begin{align*}
&\frac{\partial}{\partial t} g_\sigma - \Delta g_\sigma - 2 \, (1-\sigma) \, \Big \langle \frac{\nabla H}{H},\nabla g_\sigma \Big \rangle \\
&+ 2 \, \sum_{i=1}^n \frac{H^{\sigma-1}}{\rho+\lambda_i} \, (D_i \rho)^2 - \sigma \, |A|^2 \, (g_\sigma + K_0) \\
&\leq -H^{\sigma-1} \, \max \{\omega + |A|^2 \, \rho,0\} - \sigma \, (1-\sigma) \, H^{\sigma-3} \, (\rho - \delta \, H) \, |\nabla H|^2 + C \, H^\sigma 
\end{align*}
on the set $\{g_\sigma \geq 0\}$. Note that $g_\sigma \leq H^{\sigma-1} \, \rho$ by definition of $g_\sigma$. Since $\sigma \in (0,\frac{1}{2})$, we have $2\sigma \, \frac{g_\sigma}{\rho} \leq H^{\sigma-1}$ at each point on the hypersurface. This implies
\begin{align*}
&\frac{\partial}{\partial t} g_\sigma - \Delta g_\sigma - 2 \, (1-\sigma) \, \Big \langle \frac{\nabla H}{H},\nabla g_\sigma \Big \rangle + \sum_{i=1}^n \frac{H^{\sigma-1}}{\rho+\lambda_i} \, (D_i \rho)^2 \\
&\leq -H^{\sigma-1} \, \max \{\omega + |A|^2 \, \rho,0\} + \sigma \, |A|^2 \, (g_\sigma + K_0) + C \, H^\sigma \\
&\leq -2\sigma \, \frac{g_\sigma}{\rho} \, (\omega + |A|^2 \, \rho) + \sigma \, |A|^2 \, (g_\sigma + K_0) + C \, H^\sigma \\
&= -2\sigma \, \frac{g_\sigma}{\rho} \, \omega + \sigma \, |A|^2 \, (K_0-g_\sigma) + C \, H^\sigma
\end{align*}
on the set $\{g_\sigma \geq 0\}$. Therefore, we have
\begin{align*}
&\frac{d}{dt} \bigg ( \int_{M_t} g_{\sigma,+}^p \bigg ) \\
&\leq -p(p-1) \int_{M_t} g_{\sigma,+}^{p-2} \, |\nabla g_\sigma|^2 + 2 \, (1-\sigma) \, p \int_{M_t} g_{\sigma,+}^{p-1} \, \Big \langle \frac{\nabla H}{H},\nabla g_\sigma \Big \rangle \\
&- p \int_{M_t} \sum_{i=1}^n \frac{g_{\sigma,+}^{p-1} \, H^{\sigma-1}}{\rho+\lambda_i} \, (D_i \rho)^2 - 2\sigma \, p \int_{M_t} \frac{g_{\sigma,+}^p}{\rho} \, \omega \\
&+ \sigma \, p \int_{M_t} g_{\sigma,+}^{p-1} \, (K_0-g_\sigma) \, |A|^2 + \int_{M_t} (C \, H^\sigma \, g_{\sigma,+}^{p-1} - H^2 \, g_{\sigma,+}^p).
\end{align*}
Integration by parts gives
\begin{align*}
-\int_{M_t} \frac{g_{\sigma,+}^p}{\rho} \, \omega &\leq C \, p \int_{M_t} \frac{g_{\sigma,+}^{p-1}}{H} \, |\nabla \rho| \, |\nabla g_\sigma| \\
&+ C \int_{M_t} \frac{g_{\sigma,+}^p}{H^2} \, (|\nabla \rho|+1) \, (|\nabla H|+1) + C \int_{M_t} \frac{g_{\sigma,+}^p}{H^2} \, (|\nabla \rho|^2+1),
\end{align*}
where $C$ is a positive constant that depends only on $N$, $M_0$, $\delta$, and $T$. Putting these facts together, we obtain
\begin{align*}
&\frac{d}{dt} \bigg ( \int_{M_t} g_{\sigma,+}^p \bigg ) \\
&\leq -p(p-1) \int_{M_t} g_{\sigma,+}^{p-2} \, |\nabla g_\sigma|^2 + 2 \, (1-\sigma) \, p \int_{M_t} g_{\sigma,+}^{p-1} \, \Big \langle \frac{\nabla H}{H},\nabla g_\sigma \Big \rangle \\
&- p \int_{M_t} \sum_{i=1}^n \frac{g_{\sigma,+}^{p-1} \, H^{\sigma-1}}{\rho+\lambda_i} \, (D_i \rho)^2 + C \, \sigma \, p^2 \int_{M_t} \frac{g_{\sigma,+}^{p-1}}{H} \, |\nabla \rho| \, |\nabla g_\sigma| \\
&+ C \, \sigma \, p \int_{M_t} \frac{g_{\sigma,+}^p}{H^2} \, (|\nabla \rho|+1) \, (|\nabla H|+1) + C \, \sigma \, p \int_{M_t} \frac{g_{\sigma,+}^p}{H^2} \, (|\nabla \rho|^2+1) \\
&+ \sigma \, p \, K_0^p \int_{M_t} |A|^2 + \int_{M_t} (C \, H^\sigma \, g_{\sigma,+}^{p-1} - H^2 \, g_{\sigma,+}^p).
\end{align*}
where $C$ is a positive constant that depends only on $N$, $M_0$, $\delta$, and $T$. Using the identity
\[\frac{\nabla H}{H} = \frac{\nabla \rho - H^{1-\sigma} \, \nabla g_\sigma}{(1-\sigma) \, \rho + \sigma \, \delta \, H},\]
we obtain
\[\Big \langle \frac{\nabla H}{H},\nabla g_\sigma \Big \rangle \leq \frac{\langle \nabla \rho,\nabla g_\sigma \rangle}{(1-\sigma) \, \rho + \sigma \, \delta \, H} \leq C \, H^{-1} \, |\nabla \rho| \, |\nabla g_\sigma|\]
and
\[\frac{|\nabla H|}{H} \leq C \, \frac{|\nabla \rho|}{H} + C \, \frac{|\nabla g_\sigma|}{g_{\sigma,+}}.\]
This gives 
\begin{align*}
&\frac{d}{dt} \bigg ( \int_{M_t} g_{\sigma,+}^p \bigg ) \\
&\leq -p(p-1) \int_{M_t} g_{\sigma,+}^{p-2} \, |\nabla g_\sigma|^2 - p \int_{M_t} \sum_{i=1}^n \frac{g_{\sigma,+}^{p-1} \, H^{\sigma-1}}{\rho+\lambda_i} \, (D_i \rho)^2 \\
&+ C \, (p + \sigma \, p^2) \int_{M_t} \frac{g_{\sigma,+}^{p-1}}{H} \, |\nabla \rho| \, |\nabla g_\sigma| + C \, \sigma \, p \int_{M_t} \frac{g_{\sigma,+}^p}{H^2} \, |\nabla \rho|^2 \\ 
&+ C \, \sigma \, p \int_{M_t} \frac{g_{\sigma,+}^p}{H^2} \, (|\nabla \rho|+1) + C \, \sigma \, p \int_{M_t} \frac{g_{\sigma,+}^{p-1}}{H^2} \, |\nabla g_\sigma| \\ 
&+ \sigma \, p \, K_0^p \int_{M_t} |A|^2 + \int_{M_t} (C \, H^\sigma \, g_{\sigma,+}^{p-1} - H^2 \, g_{\sigma,+}^p).
\end{align*}
where $C$ is a positive constant that depends only on $N$, $M_0$, $\delta$, and $T$.

Consequently, there exists a positive constant $c_0$, depending only $N$, $M_0$, $\delta$, and $T$, with the following property: if $p \geq \frac{1}{c_0}$ and $\sigma \leq c_0 \, p^{-\frac{1}{2}}$, then we have
\[\frac{d}{dt} \bigg ( \int_{M_t} g_{\sigma,+}^p \bigg ) \leq C \, \sigma \, p \int_{M_t} g_{\sigma,+}^p + \sigma \, p \, K_0^p \int_{M_t} |A|^2 + \int_{M_t} (C \, H^\sigma \, g_{\sigma,+}^{p-1} - H^2 \, g_{\sigma,+}^p).\] 
Finally, since we have a lower bound for the function $H$ for $t \in [0,T)$, we obtain a pointwise upper bound for the function $C \, H^\sigma \, g_{\sigma,+}^{p-1} - H^2 \, g_{\sigma,+}^p$ for all $t \in [0,T)$. This yields 
\[\frac{d}{dt} \bigg ( \int_{M_t} g_{\sigma,+}^p \bigg ) \leq C \, \sigma \, p \int_{M_t} g_{\sigma,+}^p + \sigma \, p \, K_0^p \int_{M_t} |A|^2 + C \, p^p \int_{M_t} H^{2-(2-\sigma)p}.\] 
This completes the proof of Proposition \ref{Lp.bound.2}. \\

As above, we can now the Michael-Simon Sobolev inequality (cf. \cite{Michael-Simon}) and Stampacchia iteration to show that 
\[\rho \leq \delta \, H + C(N,M_0,\delta,T)\] 
for all $t \in [0,T)$ and all points on $M_t$. 

\section{Concluding remarks}

1. Having extended the noncollapsing estimate of Andrews \cite{Andrews} to Riemannian manifolds, we can conclude that the interior gradient estimates from Theorem 1.8' in \cite{Haslhofer-Kleiner} also hold for mean curvature flow of mean convex hypersurfaces in Riemannian manifolds. The constant in the interior gradient estimate will depend on the noncollapsing constant and also on the length of the time interval $[0,T)$.

2. In a joint work with Gerhard Huisken \cite{Brendle-Huisken}, we have recently defined a notion of mean curvature flow with surgery for mean convex surfaces in $\mathbb{R}^3$, thereby extending the work of Huisken and Sinestrari \cite{Huisken-Sinestrari3} to the case $n=2$. Using the sharp noncollapsing estimate established in this paper, it is possible to extend the surgery construction in \cite{Brendle-Huisken} to mean convex surfaces in three-dimensional Riemannian manifolds. This leads to the following result: 

\begin{theorem}[S.~Brendle, G.~Huisken]
Let $N$ be a compact Riemannian manifold of dimension $3$, and let $M_0$ be a closed, embedded surface in $N$. We assume that $M_0$ is the boundary of a domain in $N$, and has positive mean curvature. Finally, let $T>0$ be a given positive real number. Then we can define a mean curvature flow with surgery starting from $M_0$ which is defined on the time interval $[0,T)$. 
\end{theorem}

Of course, the flow with surgery can become extinct prior to time $T$.

In order to extend the analysis in \cite{Brendle-Huisken} to the case of Riemannian manifolds, we choose curvature cutoffs extremely large. This means that we only need to do surgery on necks that are extremely small. However, if the neck is sufficiently small, then the background metric on $N$ will be very close to the Euclidean metric, so that we can still apply the surgery procedure described in \cite{Brendle-Huisken}. The details will appear elsewhere.

\end{document}